\documentclass[10pt]{amsart}
\usepackage[cp1251]{inputenc}
\usepackage[english
]{babel}
\usepackage{amsmath}
\usepackage{amssymb}
\usepackage{amsfonts}
\usepackage[dvips,pdftex]{graphicx}

\newcommand{\ee}{\varepsilon}

\newcommand{\lllK}{\langle}
\newcommand{\rrrK}{\rangle}

\newcommand{\hB}{\hat{B}}
\newcommand{\hX}{\hat{X}}

\newcommand{\hx}{\hat{x}}
\newcommand{\hy}{\hat{y}}

\newcommand{\wt}{\widetilde{dist}}

\begin{document}

\title{Strongly normal cones and the midpoint locally uniform rotundity}
\author{{K. V. STOROZHUK}}
\address{Konstantin Storozhuk
\newline\hphantom{iii} Sobolev Institute of Mathematics,
\newline\hphantom{iii} acad Koptyug avenue, 4,
\newline\hphantom{iii} Novosibirsk state university,
\newline\hphantom{iii} Pirogova str, 2,
\newline\hphantom{iii} 630090, Novosibirsk, Russia}%
\email{stork@math. nsc. ru}

\maketitle


\noindent{\bf Abstract } We give the method of construction of
normal but not strongly normal positive cones. \vskip2mm



\vskip-25mm
\section{  Cones, generated by a convex set on a
hyperplane}

Let $X$ be a normed ordered space and let $K=X_+$ be its positive
cone, which orders the space $X$ in the following way: $x\leq y
\Leftrightarrow y-x\in K$. The cone generates $X$ if $X=K-K$.  An
(order) interval (or conic segment)  $\lllK a, b\rrrK$ is the set
$a+K\cap b-K$.

Let $a\in X$ and $B\subset X$. The number $dist(a,B)=\inf\{\rho
(a,b)\mid b\in B\}$ is called the distance between  $a$ and  $B$.
The number $\widetilde{dist}(A,B) =\sup_{a\in A} dist(a,B)$ is
called the nonsymmetric distance from the set $A$ to the set $B$.
The distance between  $A$ and $B$ is defined to be the number
$dist(A,B)=\min\{\widetilde{dist}(A,B),\widetilde{dist}(B,A)\}. $

Krein \cite{Krein} introduced the notion of a normal cone.
In our terminology, a cone is normal if the function
$\rho(x,y)=dist (\lllK 0,x\rrrK,\lllK 0,y \rrrK)$, defined on the
set $K\times K$, is continuous at $(0,0)$.

In \cite{EmWo}, there appeared the strong normality condition: it
means the continuity of the function $\rho$ on the entire $K\times
K$. It was also indicated there that it was not known whether every
normal cone is strongly normal. See also \cite{Em} where the strong
normality condition is actively used in Chapter 2.1.

It is easy to see that the function $\rho$ is continuous on $K\times K$
if it is continuous at the points of the form $(x,x)$ with
respect to one argument, i.e., when for any $x\in K$ $\rho(x,
y_n)\to 0$ for $y_n\to x$. In what follows, the phrase like
"function $\rho$ is continuous at the point $x$"\ will mean exactly
that.

It is also easy to see that, in checking continuity of the function
$\rho$, it suffices to restrict ourselves to the analysis of its
continuity at the points of a hyperplane section $\hB$, i.e., to
consider the sequences of the form $\hy_n\in \hB,$ $\hy_n \to \hx$.

In the language of multivalued maps, strong normality is the
continuity of the map
 $z\mapsto \lllK 0,z\rrrK$ at all the points  $z\in K$.
One can define a condition of "semi-strong"\ normality,
which involves semicontinuity of the function
$z\mapsto \lllK 0, z\rrrK$. Semicontinuity from above means that for
$\hy_n\to \hx$ $\wt (\lllK 0,\hy_n\rrrK, \lllK 0, \hx\rrrK)\to 0$,
whereas semicontinuity from below, that  $\wt (\lllK 0,\hx\rrrK,
\lllK0, \hy_n\rrrK)\to 0$.

In the present paper we give some examples of non-strongly normal
cones.
Our cones will be of the form $K=\{t\hx\mid \hx\in \hB,\ t\geq 0\}$,
where $\hB$ is a convex set in the hyperplane section of $K$.
These are good cones. We characterize strong
normality of the cone $K$ in terms of geometry of the generating
section $\hB$.

Throughout this paper, $X$ is a real Banach space, $B\subset X$ is a
bounded convex closed set in $X$, $E=\Bbb R\times X$, $K\subset \Bbb
E$ is the cone generated by the set $\hB:=\{1\}\times B$. It is
clear that such cone is normal and if its interior is non-empty then
$K$ generates the space $E$. We define a norm in $E$ by: $\|(t\times
x)\|=|t|+\|x\|$. If $x\in X$ then by $\hat{x}$ we will denote the
vector $(1\times x)\in E$.

A segment $[a,b]$ in a vector space will be understood as the
ordinary algebraic segment $\{a+t(b-a)\mid 0\leq t\leq 1\}$. In
particular, $[0,b]=\{tb\mid t\in [0,1]\}$.

Clearly, if  $z\in K$ then $[0,z]\subset \lllK 0,z\rrrK$. By the
thickness of the interval  $\lllK 0, z\rrrK$ we will mean the
distance from $\lllK 0, z\rrrK$ to the segment $[0,z]\subset \lllK
0,z\rrrK$.

It is not difficult to observe that  $x$ is an extreme point of $B$
if and only if the thickness of the interval  $\lllK 0,\hat{x}\rrrK$
equals zero, i.e., $[0,\hat{x}]= \lllK 0,\hat{x}\rrrK$ (see, for
example, \cite[Definition 1. 42 and Lemma 1. 43.]{ACT}).

\bf Theorem  1\it. Let $x\in B$. If the point $x$ is not extreme but
there exist extreme points $y_n\to x$ then the function $\rho$ is
discontinuous at $\hat{x}$.

Proof\rm. The intervals $\lllK 0,\hat{y}_n\rrrK$ coincide with the
segments $[0,\hat{y}_n]$ and converge in the limit to the segment
$[0,\hat{x}]$, while the interval
 $\lllK 0,\hat{x}\rrrK$ is "thick"\ , i.e., it contains, besides the segment
$[0,\hat{x}]$, some extraneous points. The rest is obvious.

It is already in $\Bbb R^3$ where there exists a convex closed set
$B$ whose extreme points do not form a closed set\rm. The reader
can construct such an example him/herself or find an example, say,
in \cite{Rockaf}, Chapter~4, the example after Corollary 18.5.3.

\bf Corollary\it. In $\Bbb R^4$ there exists a normal but not
strongly normal cone\rm.

It is intuitively clear that all the unexpected «metamorphoses» of
the intervals must take place on the boundary of the cone; the
following theorem confirms this. The proof is not hard, but rather
long. It is given at the end of the article.

\bf Theorem 2\it. Let $B$ be a bounded convex subset of $X$. The
function $\rho(x,y)$ is continuous at the interior points of the
cone $K\subset E$\rm.

\bf Lemma on an extraneous point\it. Let  $x\in B$ and $a\in X$.

1. The following conditions are equivalent:

(a)  $[x-a, x+a]\subset B$,

(b) $\frac{\hat{x}+a}{2}\in\lllK 0,\hat{x}\rrrK\subset E$.

2. If $l(x)$ is the supremum of the lengths of the segments
$B$ with the midpoint at $x\in B$ then the thickness of the interval
$\lllK 0,\hat{x}\rrrK$ is no less than
$\frac{\|l\|}{4}$ and no greater than  $\|l\|$\rm.

The proofs are rather simple and omitted here; we only present two
illustrations. \vskip1mm

\hskip-5mm\includegraphics[scale=0.35]{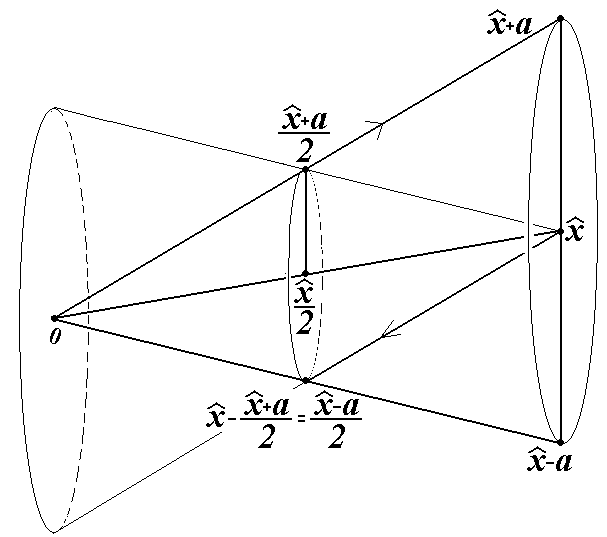} \vskip-50mm
\hskip68mm
\includegraphics[scale=0.35]{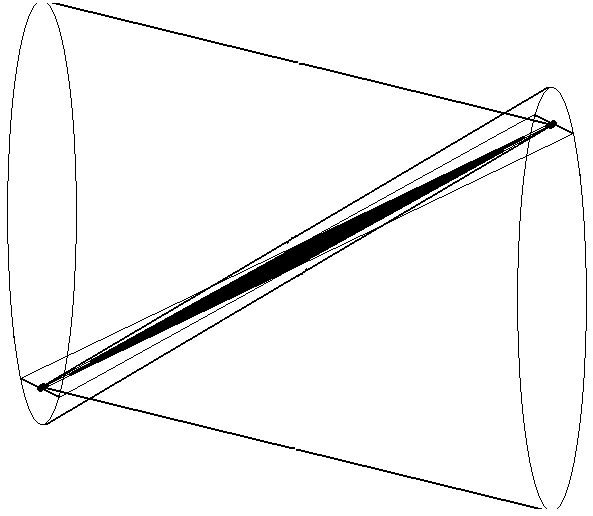}

\vskip12mm\bf Lemma on the midpoints of  long chords\it. Let $x\in
B$ be an extreme point of the set  $B$. The following are
equivalent:

(a) there exists a sequence converging to $x$ and consisting of ``uniformly
non-extreme'' points  $y_n\in B$, i.e., the
midpoints of the chords $B$ of the lengths greater than some $l>0$;

(b) the function $\rho$ is discontinuous at $\hat{x}$.

Proof\rm. Since  $x$ is an extreme point of $B$, we have $\lllK
0,\hx\rrrK=[0,x]$.

Let $\hy_n\to \hx$. The segments $[0,\hy_n]$ converge to the
segment $[0,\hx]=\lllK 0,\hx\rrrK$. Therefore, the distance from the
intervals $\lllK0,\hy_n\rrrK$ to  $\lllK 0,\hx\rrrK=[0,\hx]$
converges to zero if and only if the thickness of the intervals
$\lllK0,\hy_n\rrrK$ tends to zero. It remains to recall the second
part of  Lemma on an extraneous point. The lemma is proved.

\bf Theorem on semicontinuity\it.

1. If $\dim X<\infty$, then the function $\rho$ is semicontinuous
from below regardless of the
 set $B$  generating the cone; i.e., the limit interval can only ``grow up''.

2. If $X$ is strongly convex, i.e., all the points of the unit
sphere $S$ are extreme points of the ball $B$,
 then the function $\rho$ is semicontinuous from above, i.e., the limit interval can only
decrease\rm.

The first part of the theorem can be derived 
from the compactness of the finite-dimensional ball and the
reasoning that the reader carried out while proving the Lemma on an
extraneous point. The second part follows from the fact that the
intervals $\lllK 0,\hx\rrrK$  coincide with the segments $[0,\hx]$
for $x\in S$ and so on; cf. with the proof of Lemma on the long
chords.\vskip-2mm

\section{ Strong normality and the midpoint locally
uniform rotundity }

A space is called uniformly convex if $x_n\in S, y_n\in S$  and
$\|\frac{x_n+y_n}{2}\|\to 1$ imply  $\|x_n-y_n\|\to 0$. A space is
called locally uniformly convex \cite{Lova} if, in the previous
definition, we additionally ``fix'' one of the endpoints of the chords:
$x_n\equiv x\in S$.

If, on the other hand, we fix not an endpoint of the chords but
their midpoints then we arrive at what K.W. Anderson  \cite{An}
called the \it midpoint locally uniform rotundity \rm (MLUR).
(Anderson or somebody else attributed this property to G. Lumer and
M. M. Day).

Namely,
$X$ is midpoint locally uniformly rotund
(MLUR) if whenever $x_n\to x\in S,$ $\|x\pm v_n\|\to 1$, we have
$v_n\to 0$. It means geometrically that any point  $x$ of the sphere
$S$ is uniformly far from the midpoints of  long chords of this
sphere. It is precisely this property that is used in  Lemma on the
midpoints of  long chords. It follows from this Lemma that, at the
extreme point $x$ of the ball, the MLUR condition is equivalent to
the continuity of the function $\rho$ at the point $\hx$. Since at
the interior points of the cone the function $\rho$ is always
continuous (Theorem 2), the following theorem holds:

\bf Theorem 4\it. Let $X$ be strongly convex and let $B$ be the unit
ball in $X$. A cone $K\subset E$ is strongly normal if and only if
$X\in MLUR$. \rm

A substantial number of papers are devoted to the research on the MLUR property, its comparing
with other characteristics of convexity of the sphere, the duality
issues and the possibilities of the MLUR and non-MLUR
renormalization. It is known that any separable Banach space is
isomorphic to a locally uniformly convex and, therefore, MLUR space.
On the other hand, the majority of ``decent''\ spaces admit a non-MLUR
renormalization. For example, the separable spaces, containing
$l_1$, possess a strongly convex not-MLUR norm. See, \cite{Kadec1},
\cite{Kadec2}. Let us refer, for example, to the paper \cite{Smit};
there one can find three examples of non-MLUR spaces. Let us give
our own example for the sake of completeness of the exposition.

\bf Example. \rm  A convex closed subset  $B$ in a Hilbert space
$H$ with an extreme point which is approximated by the centers of
the infinite dimensional subdiscs.

Let $e_i$, $i=1,2,\ldots$ be a Hilbert basis in  $H$. Let $B_n$ be
the unit ball in the subspace spanned by the  $e_{n+1},
e_{n+2},\ldots$. For $n\geq 1$ \  we set $Z_n=\{te_1+B_n
\mid|t|\leq\frac{n}{n+1}\}$. Let $B$ be the closure of the convex
hull of all $Z_n$. The interior of $B$ is non-empty (since the
interior of $Z_1$ is non-empty: it contains a ball in $H$ of the
radius $\frac12$). The sequence $y_n=\frac{n}{n+1}e_1$ of the
centers of the (infinite dimensional!) discs $y_n+B_n\subset B$
converges to the exterior point $e_1$ of the set $B$. Being
considered as the unit ball, the set $B$ defines in $H$ an
equivalent non-MLUR norm.

\vskip-2mm\section{ Appendix. Proof of Theorem~2}

Let us note that if  $K_1\subset K_2$ are two cones then for all $x$
$\lllK 0,x\rrrK_{K_1}\subset\lllK 0,x\rrrK_{K_2}$. Secondly, if
$K'=T(K)$ where $T$ is a linear transformation then $T\lllK
0,x\rrrK_{K}=\lllK 0,Tx\rrrK_{K'}$.

Now, let  $\hat{x}\in int(K)$ and $\hat{y}_n\to \hx$.  Without any
loss of generality, we can take $\hx=(1\times 0)\in\hB$ (however, we
keep in mind that $B$ is not necessarily a ball). We are going to
prove that $\lllK 0,\hat{y}_n\rrrK \to \lllK 0, \hx\rrrK$.

Let $\ee>0$. Consider linear transformations $T_{+\ee}$ and
$T_{-\ee}$ of the space $\Bbb R\times X$  which preserve the first
coordinate and stretch the ``second coordinate'', i.e., the hyperplane
$X$, by $1\pm\ee$ times: $f_{\pm\ee}(t\times
z)=(t\times(1\pm\ee)z)$. The cone $K$ turns out to be in between the
cones $f_{\pm\ee}(K)$.

The sets $B_\pm=(1\pm\ee)B$ are separated from the set  $B\subset
X$; moreover, the distance from the boundaries of the sets  $B_\pm$
to the boundary of the set $B$ is no less than the number
$\ee\cdot\delta$, where  $\delta$ is the radius of the ball with the
center at the point  $x=0$, which is contained in $B$. Let us show
this in the picture.

\hskip-\parindent
\parbox[t]{189pt} { Let $l$ be a
supporting plane at the point $D\in B$, $D^+=(1+\ee)D$. Clearly, the
distance from $D^+$ to the set $B$ is no less than $|CD^+|$. It
follows from the similarity of the triangles $DAO$ and $DCD^+$ that
$\frac{CD'}{OA}=\frac{DD^+}{OD}=\ee$, therefore, $|CD^+|=\ee|OA|\geq
\ee\delta$. A similar argument shows that the distance from  $D$ to
$B_-$ is also is at least $\ee\delta$. Thus, the set $B$ is included
in the interval $B_-\subset B\subset B_+$ together with its
$\ee\delta$-neighborhood. }

\vskip-41mm\hskip190pt\includegraphics [scale=0.18]{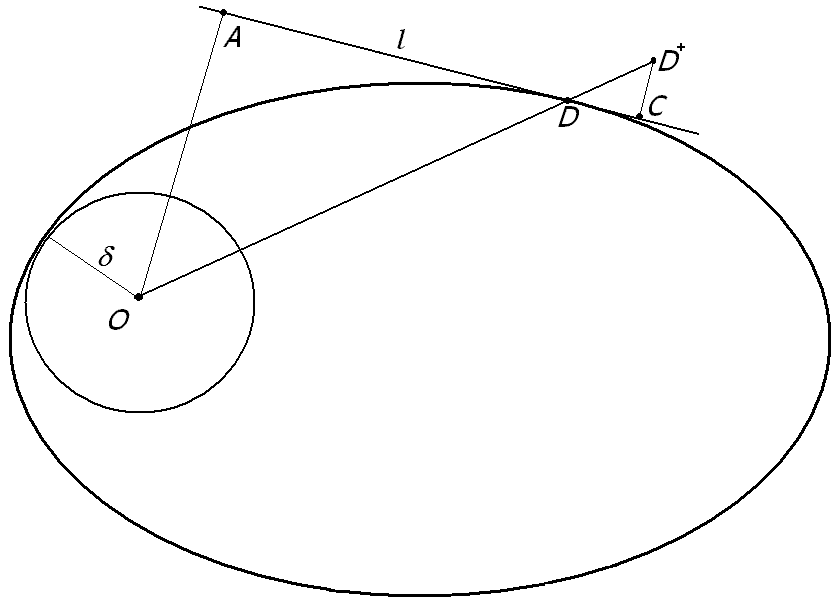} \vskip3mm

Let us consider further the linear transformations $G_n:\Bbb R\times
X\to \Bbb R\times X$ which also preserve the first coordinate and
move affinely the hyperplane
 $\hX$ in such a way that the point
$\hat{y}_n$ goes to $0$. Formally, $G_n(t\times z)=(t\times
(z-ty_n))$. Denote order intervals defined by the cones $G_n(K)$ by
the symbol $\lllK ,\rrrK_n$; whereas, the ones defined by the cones
$f_{\pm\ee}(K)$, by the symbols $\lllK ,\rrrK_{\pm\ee}$. It follows
from the previous argument that if $\|y_n-x\|<\ee\delta,$ then
$K_{-\ee}\subset K_n\subset K_{+\ee}$. Therefore, for the intervals
we  have
$$T_{-\ee}(\lllK 0,x\rrrK)=\lllK 0,x\rrrK_{{-\ee}}\subset \lllK 0,x\rrrK_n \subset \lllK
0, x\rrrK_{{+\ee}}=T_{+\ee}(\lllK 0,x\rrrK). $$ Recall now that
$\lllK 0,x\rrrK_n=G_n\lllK 0,G_n^{-1}x\rrrK=G_n\lllK 0,y_n\rrrK$.
However, for large $n$ the intervals $G_n(\lllK 0,y_n\rrrK)$ are
close to the intervals  $\lllK 0,y_n \rrrK$. Thus, for large $n$
$T_{-\ee}(\lllK 0,x\rrrK)\subset \lllK 0,y_n\rrrK \subset
T_{+\ee}(\lllK 0,x\rrrK)$. Since, for small $\ee$, the intervals
$f_\pm \lllK 0, x\rrrK$ are close to the interval $\lllK 0,x\rrrK$,
the theorem is proved.


\begin{thebibliography}{}

\bibitem{Krein}Krein, M. Proprie'te's fondamentales des ensembles coniques normaux
dans l'espace de Banach. (French) C. R. (Doklady) Acad. Sci. URSS
(N. S. ) 28, (1940). 13–17.

\bibitem{EmWo}Emelyanov, E. Yu. ; Wolff, M. P. H. Positive operators on Banach
spaces ordered by strongly normal cones.  Positivity 7 (2003), N.
1-2, 3–22.

\bibitem{Em}  Emel'yanov. E.Yu. Non-spectral asymtotic analysis of
one-parameter operator semigroups. Operator theory Advances and
applications, vol.173. Birkhauser 2007.

\bibitem{ACT}Aliprantis, Charalambos D. ; Tourky, Rabee.
Cones and duality. Graduate Studies in Mathematics, 84. American
Mathematical Society, Providence, RI, 2007.


\bibitem{Rockaf} 
R.T. Rockafellar. Convex Analysis. Princeton University Press,
Princeton,. N.J., 1970.

\bibitem{Lova} A. R. Lovaglia. Locally uniformly convex Banach spaces, Trans.
Amer. Math. Soc. 78, (1955), 225-238.

\bibitem{An} K. W. Anderson. Midpoint local uniform convexity, and other geometric properties
of Banach spaces, Ph. D. dissertation, Univ. Illinois, Urbana, IL,
1960.


\bibitem{Kadec1} Kadec\!', M. \v{I}. Spaces isomorphic to a locally uniformly convex space. (Russian)
Izv. Vys\v{s}. U\v{c}ebn. Zaved. Matematika 1959 1959 no. 6 (13),
51–57.

\bibitem{Kadec2}Kadec\!', M. \v{I}. Relation between some properties of convexity of the unit ball of a
Banach space. Funct. Anal. Appl. 16, N 3, 204-206.

\bibitem{Smit} Smith, Mark A. Some examples concerning rotundity in
Banach spaces. Math. Ann. 233 (1978), no. 2, 155–161.


\end{thebibliography}
\end{document}